\documentclass[a4paper, 12pt, oneside]{amsart}

\usepackage{amsmath,amssymb,citesort,enumerate,times,latexsym}
\usepackage{amsthm}
\usepackage[all]{xy}

\def\dfn{\emph}                
\def\ie{i.e.}                  

\newcommand{\kitty}[1]               
    {\ensuremath{\mathcal{#1}}}      
\newcommand{\cat}[1]                 
    {\textbf{#1}}

\def\Ker{\mathrm{Ker}\,}             
\def\isom{\ensuremath{\cong}}        
\def\comp{\ensuremath{\circ}}        

\def\pt{\ensuremath{\partial}}       

\newcommand{\act}[2]{\ensuremath{{ }^{#1}\!#2}}

\newcommand{\oset}[1]{\langle#1\rangle}

\def\zee{\mathbb{Z}}

\title{Group Objects and Internal Categories}
\author{Magnus Forrester-Barker}
\email{map601@bangor.ac.uk}
\address{Ysgol Gwybodeg, Prifysgol Cymru Bangor, Bangor, Gwynedd, UK, LL57 1UT}
\date{4 December 2002}

\begin{document}

\begin{abstract}
Algebraic structures such as monoids, groups, and categories
can be formulated within a category using commutative diagrams.  In
many common categories these reduce to familiar cases.  In particular, group
objects in \cat{Grp} are abelian groups, while internal categories in
\cat{Grp} are equivalent both to group objects in \cat{Cat} and to crossed
modules of groups.  In this exposition we give an elementary introduction to
some of the key concepts in this area.
\end{abstract}

\maketitle

This expository essay was written in the winter of 1999-2000, early in the
course of my PhD research, and has since been updated with supplementary
references.  I hope you will find it useful.  I am indebted to my supervisor,
Professor Tim Porter, for his help in preparing this article.

\section{Groups Within a Category}

Let \kitty{C} be a category with finite products.  For this it is
necessary and sufficient that \kitty{C} have pairwise products
(i.e. for any 2 objects $C,D \in Ob(\kitty{C})$, there is a product $C \times D$)
and a terminal object, which we shall denote by $\mathbf{1}$.  Examples
of suitable categories include $\cat{Set}, \cat{Grp}, \cat{Top}$ and
$\cat{Ab}$.

Let $G$ be an object of \kitty{C}.  Then $G\times G$ is also an object of
\kitty{C}.  Suppose we can find a morphism $m:G\times G \rightarrow G$ such
that the diagram 
\[
\xymatrix{
  G\times G\times G \ar[r]^-{id_G\times m} \ar[d]_-{m\times id_G} & G\times G 
\ar[d]^-m \\
  G \times G \ar[r]_-m & G
}
\]
commutes.  Now, $m$ is a binary operation and if we temporarily think in
terms of elements we see that the diagram ensures that the operation is 
associative.  For example, take $\kitty{C}=\cat{Set}$; then $G$ is a set 
and $(id_G\times m)(a,b,c):= (a,bc)$ for $a,b,c \in G$, etc.  We may
take this diagram as a general definition of associativity, even for
categories in which the objects do not have elements.  Therefore $\oset{G,m}$
gives us an abstract semigroup in \kitty{C}.  It is often convenient, in 
categories such as \cat{Set}, to think in terms of the elements and their 
images under the given morphisms\footnote{Sometimes it is actually more 
confusing to think in terms of elements, for example when dealing with 
opposite categories.}, but in fact the definitions are much more general.

Suppose we also have a morphism $e:\mathbf{1} \rightarrow G$ from the terminal
object to $G$, such that 
\[
\xymatrix{
  \mathbf{1}\times G \ar[r]^-{e\times id_G} \ar[dr]_-{\isom} & G\times G \ar[d]^m 
  &  G\times \mathbf{1} \ar[l]_-{id_G\times e} \ar[dl]^-{\isom} \\
      & G &
}
\]  
commutes\footnote{Note that the isomorphisms $\mathbf{1}\times G \isom G \isom
G\times \mathbf{1}$ follow automatically from the definition of the product
and the fact that $\mathbf{1}$ is a terminal object  (so there is a unique
arrow $G \rightarrow \mathbf{1}$).}.  Then $e$ acts as an identity for $G$.
For example, in $\cat{Set}$, the terminal object $\mathbf{1}$ is a singleton 
and so $e$ picks out an element of $G$ which acts as the identity under the
``multiplication'' $m$.  

Hence an object $G$ with morphisms $m, e$ as above can be thought of as a
monoid in \kitty{C}.  In fact we can define (see MacLane \cite[section 3.6]
{MacLaneS:catwm}) a monoid in \kitty{C} to be precisely a 
triple $\oset{G, m, e}$ such that the above diagrams commute.  In the case 
of $\cat{Set}$, the monoids thus given correspond exactly to monoids in the 
classical, i.e. non-categorical, sense.  It is clear that, given any classical 
monoid (written multiplicatively), we can define $m$ by $m(x,y):=xy$ and if 
$1$ denotes the identity element, we define $e:\star
\longmapsto 1$ (taking $\mathbf{1}=\{\star\}$).  Conversely, given any 
categorical monoid in $\cat{Set}$, the morphism $m$ gives a well-defined 
associative binary operation for which $e(\star)$ is an identity.

To get from here to groups it remains to define inverses.  Again we must do
this using morphisms, as the objects in our category need not have
elements.  Let $i$ be a morphism $G \rightarrow G$ such that\[
\xymatrix{ G \ar[r]^-{\Delta} \ar[d]_{\exists !} & G\times G \ar[r]^{i \times
   id_G} & G\times G \ar[d]^m \\
\mathbf{1} \ar[rr]_e & & G
}
\]
commutes, where $\Delta: G \rightarrow G\times G$ is the diagonal morphism
(i.e. $p\Delta = q\Delta = id_G$, where $p,q$ are the projections from the
product to its components).  In $\cat{Set}$, if $g\in G$ then 
$\Delta (g)=(g,g)$ and the commutativity of the diagram gives 
$m(i(g),g)=e(\star)$, so $i(g)$ is a left inverse for $g$.  We can also 
replace $i\times id_G$ in the diagram with $id_G \times i$, which would 
give $i(g)$ as a right inverse for $g$.  Thus to get a two-sided inverse 
for each element in $G$, we may stipulate that $i$ should exist so that both 
diagrams commute.

We are now in a position to define a \dfn{group object} in \kitty{C} to be
an ordered quadruple $\oset{G, m, e, i}$ with morphisms defined as above.
In $\cat{Set}$, the group objects are just groups in the classical sense,
as we have seen.  In $\cat{Top}$, each set $G$ is a topological space and
all maps are continuous, so the group objects are topological groups.  In the
category of differentiable manifolds\footnote{The morphisms are smooth maps of 
class  $C^\infty$, \ie they are infinitely differentiable.},
the group objects are Lie groups.  In $\cat{Grp}$, the objects
are themselves groups; in this case, the group objects are abelian groups.
Why?

Firstly, every abelian group is a group object in $\cat{Grp}$.  Suppose $A$
is an abelian group.  Then it consists of a set $A$ together with an
associative, commutative binary operation with an identity and inverses.  We
can view the binary operation as a map $m:A\times A \rightarrow A$ and the
identity can be considered as a map $e:\mathbf{1}\rightarrow A$ which selects
the identity element; similarly there is a map $i:A\rightarrow A$ which takes 
each element in $A$ to its inverse.  These maps then clearly satisfy the 
diagrams given above, so that $\oset{A, m, e, i}$ is a group object in 
$\cat{Set}$.  However, $A$ is a group, so it is an object of $\cat{Grp}$.  
It remains to show that $m, e, i$ are group homomorphisms and hence morphisms in
$\cat{Grp}$.  Now $A\times A$ is a direct product of groups so its
multiplication is defined componentwise.  Writing $\otimes$ for the
multiplication in $A\times A$ we thus have
\[
(x,y)\otimes(x', y') := (m(x,x'), m(y,y')).
\]
We now get 
\[
\begin{array}{rclcl}
m((x,y)\otimes(x',y')) & = & m(m(x,x'),m(y,y')) & = &(x+x')+(y+y') \\
                       & = & (x+y)+(x'+y')       & = & m(m(x,y),m(x',y')).
\end{array}
\]
Note that here we switched notation, writing $m(x,y)=x+y$, so that the
associativity and commutativity could be visualised more easily.  Comparing
the left and right hand ends of this chain of equalities shows that $m$ is
a homomorphism, as required.  The other two are, fortunately, easier to write
down.  We have $e:\{\star\}\rightarrow A$; now $\{\star\}$ is the trivial
group.  Thus we get $e(\star\star) = e(\star) = id_A = m(id_A, id_A) = 
m(e(\star), e(\star))$ and hence $e$ is a homomorphism from the trivial
group.  Finally, $i:A\rightarrow A$.  We have (again writing additively for 
convenience) $i(a+b)=-(a+b)=(-a)+(-b) =i(a)+i(b)$, so this too is a 
homomorphism.  Since these maps are all homomorphisms, they exist in 
$\cat{Grp}$ and hence we get that $\oset{A, m, e, i}$ is a group object 
in $\cat{Grp}$, as required.

Conversely, suppose we have a group object $\oset{G, m, e, i}$ in $\cat{Grp}$.
We must show that $G$ is an abelian group.  Now $G$, being a group, already
has an operation which we shall write as $x\cdot y$.  This is associative, with
inverses and an identity, but not necessarily commutative.  We have a second
operation defined by the homomorphism $m$, which we shall write as
$x\ast y:= m(x,y)$.  $G$ has an identity, denoted\footnote{Note that we are
using $1_G$ to denote the identity element in the group $G$ and $id_G$ to 
denote the identity homomorphism $G\rightarrow G$.  Similar notation is
applied throughout.}
$1_G$, under $\cdot$\,; it also has an identity under $\ast$ given by 
$e:\mathbf{1}\rightarrow G$.  Since $e$ is a homomorphism and $\star$ is 
the identity element of the group $\{\star\}$, we deduce that $e(\star)=1_G$.
Now, for $x,y,z,w \in G$, we have $m(x,y)\cdot m(z,w) = 
m((x,y)\otimes (z,w))$, since $m$ is a homomorphism.  But, by definition 
of the product $\otimes$ in $G\times G$, this is just $m(xz,yw)$.  Hence, 
switching notation, we get the \dfn{interchange law}:\label{inter}
\[ 
      (x\ast y)\cdot (z\ast w) = (x\cdot z)\ast (y\cdot w). 
\]
Using this, we get:
\[
\begin{array}{rclclcl}
    x\cdot y &=& (x\ast 1_G)\cdot (1_G\ast y) &=& (x\cdot 1_G)\ast (1_G\cdot y)
    &=& x \ast y \\
    &=& (1_G\cdot x)\ast (y\cdot 1_G) &=& (1_G \ast y)\cdot (x\ast 1_G)
    &=& y \cdot x.
\end{array}
\]
Hence $\oset{G,\cdot}$ is abelian, as required.  Note also that the first line
shows that the operation given by $m$ is necessarily the same as the group
operation already defined on $G$.  

\subsection*{Group Actions}
Having formulated groups within \cat{Set}, we can also formulate group
actions by means of suitable diagrams.  Let $\oset{G,m,e,i}$ be a group 
and $X$ a set.  Then a (left) action of $G$ on $X$ is given by a map
$n:G\times X \rightarrow X$ such that the following diagrams commute:
\[
\xymatrix{ G\times G\times X \ar[r]^-{id_G\times n} \ar[d]_{m\times id_X} &
G\times X \ar[d]^n \\
G\times X \ar[r]_-n & X }
\]\[
\xymatrix{\mathbf{1}\times X \ar[r]^-{e\times id_X} \ar[dr]_-{\isom} &
G\times X \ar[d]^n \\
  & X }
\]
In the usual notation for (left) group actions, the first of these diagrams
gives $\act{g}{(\act{g'}{x})} = \act{(gg')}{x}$, while the second gives 
$\act{1_G}{x} = x$.  Hence the action so defined is just a group action in 
the usual sense and indeed the two sorts of action are entirely equivalent.  

A similar construction can be made for group objects in other categories.  For
instance in \cat{Top}, we get continuous actions of topological groups on
spaces.  The definition of an action did not make use of the group inverses,
so we can talk about more general monoid actions using the same definition.
One case worth mentioning is in the category \cat{Ab}, where a monoid
$\oset{M, m,e}$ relative to the tensor product, $\otimes$, and the ``terminal
object'' $\zee\ $ gives us a ring\footnote{Note that this is not the monoid
object defined above, since $\otimes$ is not the product (defined as a limit)
in \cat{Ab}, nor is $\zee\ $ a terminal object.  However
$\oset{Ab,\otimes,\zee.}$ forms a monoidal category, so the given diagrams
work with $\times, \mathbf{1}$ replaced by $\otimes, \zee$ respectively.  See
MacLane \cite [introduction, ch. 7]{MacLaneS:catwm}}.  In this case the action
of $M$ on an abelian group $A$ makes $A$ into a left $M$-module.

Right actions of groups (and monoids) can be defined similarly with a map
$X\times G \rightarrow X$.  These are not intrinsically different from left
actions, but the notation of right actions (for example, $x^g$) is sometimes
more natural.

\section{Other Structures}
In the same way that we have defined groups, we could define other structures
such as rings and lattices within a category by giving suitable morphisms and
commutative diagrams.  This can, in principle, be done for any algebraic 
structure that is defined in terms of operations satisfying specified
properties expressed by equations.

Also, we can dualize the construction of group objects by reversing all the
arrows in the given diagrams and replacing products by coproducts and terminal
objects by initial objects.  This gives us an operation $w:C\rightarrow C\amalg C$ such that the diagram

\[
\xymatrix{
      C \ar[r]^-w \ar[d]_-w & C \amalg C \ar[d]^{w \amalg id_C} \\
      C\amalg C \ar[r]_-{id_C \amalg w} & C \amalg C \amalg C
}
\]
commutes (giving ``co-associativity'').  Similar analogues exist to the other
diagrams given earlier, noting that we now have $\eta:C\rightarrow I$
(where $I$ is the initial object of \kitty{C}).  In place of the diagonal
morphism in the inverses diagram, we use a morphism $\nabla: C\amalg C
\rightarrow C$ which maps both copies of an element in the coproduct to the
same element in $C$.  This definition gives us a \dfn{cogroup object} in
\kitty{C} (see Rotman \cite[chapter 11]{RotmanJ:intat} for more details).  In a similar
way we could dualize any of the other algebraic structures within a category,
to get comonoids, corings, colattices, etc.  

These dual constructions do not seem to correspond to any standard algebraic
structures with which I am familiar (although presumably a cogroup in 
$\cat{Set}^{op}$ would correspond to a group in $\cat{Set}$, etc.).  They 
do find uses in modern mathematics, however.  A quick search of the MathSciNet 
database\footnote{\textit{http://klymene.mpim-bonn.mpg.de/mathscinet/} } 
reveals a number of papers on cogroups and corings (although only one 
mentioning colattices).  Putting together coalgebra and algebra structures, 
along with some further conditions yields Hopf Algebras, which in turn
can be extended to give quantum groups \cite{MajidS:fouqgt}.

\section{Internal Categories}
Just as groups and other algebraic structures can be defined in a category
by giving suitable objects and morphisms, we can sometimes build categories
within a category.  Let \kitty{C} be a finitely complete category (in fact, it
is sufficient for \kitty{C} to have pullbacks).  Suppose there are objects
$A, O \in Ob(\kitty{C})$ and morphisms
\[
\xymatrix{ A \ar@<1ex>[r]^s \ar [r]_t & O \ar@<1ex>@/^2ex/[l]^e }
\]
such that $se = id_O = te$.
We can consider $A$ as a collection of directed edges, $O$ as a collection
of vertices and the maps $s,t,e$ as giving respectively the source and target 
vertices of each edge and a loop at each vertex.  Thus $\oset{A,O,s,t,e}$ 
specifies an internal reflexive directed graph (or ``digraph'') in \kitty{C}.

Recall that every (small) category has an underlying digraph and conversely a
reflexive digraph can be turned into a category by considering vertices as
objects and edges as morphisms and specifying a composition (this may 
necessitate adding further edges).  The loops given by the reflexive property
become identity morphisms on each object.  Similarly we now seek to extend
our internal digraph in \kitty{C} to an internal category by defining a 
suitable composition.

We can form the pullback square:
\[
\xymatrix{ A\, {{}_t\times_s} A \ar[r]^-{p_2} \ar[d]_{p_1} & A \ar[d]^s \\
           A \ar[r]_t                           & O }
\]
The pullback object $A\, {{}_t\times_s} A$ can be considered as the collection
of all composable pairs of morphisms.  In a category such as \cat{Set} or
\cat{Grp} where the objects have elements and the product is an ordered
set of elements (possibly with some algebraic structure imposed) we get
$A\, {{}_t\times_s} A = \{ (f,g)\in A\times A : tf = sg \}$.  In order to 
form a category (with objects $O$ and morphisms $A$) inside \kitty{C}, we 
need to define a composition $m:A\, {{}_t\times_s} A \rightarrow A$ which is
associative and respects identities; note in particular that $m$ is also a
morphism in \kitty{C}.

I shall formulate internal categories for a category \kitty{C} in which the
objects contain elements we can work with, although in principle they could be
formulated for more general categories.  In particular, we can describe
pullbacks in terms of their elements, as above.  Define the morphism $m$ as in
the previous paragraph so that $t(m(f,g)):=tg,\; s(m(f,g)):=sf$.  It will be
convenient to write $g\comp f$ for $m(f,g)$.  Now, using the maps $tm:A\,
{{}_t\times_s} A \rightarrow O$ and $s:A\rightarrow O$, we can form the
pullback object $(A\, {{}_t\times_s} A)\, {{}_t\times_s} A$, which in terms of
elements is the set $\{ (f,g,h)\in A\times A\times A : tf = sg, t(g\comp
f)=sh\}$.  Similarly we can form the pullback object $A\, {{}_t\times_s} (A\,
{{}_t\times_s} A) = \{ (f,g,h)\in A\times A\times A : tf = s(h\comp g), tg
=sh\}$.  Since $s(h\comp g)=sg$ and $t(g\comp f)=tg$ we deduce that these two
pullback objects are in fact equal.  Hence we can form the diagram
\[
\xymatrix{ 
  A\, {{}_t\times_s} A\, {{}_t\times_s} A \ar[r]^-{m\times id_A} \ar[d]_{id_A
    \times m} &      A\, {{}_t\times_s} A \ar[d]^m \\
  A\, {{}_t\times_s} A \ar[r]_m & A 
}
\]
To get associativity of composition we now just require that this diagram
commute. 

Our final requirement for a category is that the composition respect
identities.  The morphism $e:O\rightarrow A$ selects the identity arrow
$id_x$ for each object $x$ in $O$.  If $f\in A$ with $sf=x,\,tf=y$ then
we need $id_y\comp f = f = f\comp id_x$.  Now, using the composite map 
$te = id_O$ we can form the pullback $O \,{{}_{id_O}\!\times_s} A = 
\{(x,f) : x = sf \}$.  Similarly with $se = id_O$ we get 
$A {{}_t\times_{id_O}} O = \{(f,y) : tf = y \}$.  These pullbacks have 
obvious projections onto A, namely $p:O \,{{}_{id_O}\!\times_s} A\rightarrow 
A$ with $p(x,f)=f$ and $q:A {{}_t\times_{id_O}} O\rightarrow A$ with 
$q(f,y)=f$.  These are clearly bijective, since we can define inverses 
$p^{-1}(f):=(sf,f),\ q^{-1}(f):=(f,tf)$.  Then $pp^{-1}(f)=f=id_A(f)$ and
$p^{-1}p(x,f)=(sf,f)=(x,f)=id_{O \,{{}_{id_O}\!\times_s} A}$, so $p$ is an
isomorphism; similarly for $q$.  Using these maps we express our requirement 
for identities in the commutativity of the following diagram:
\[
\xymatrix{ O\,{{}_{id_O}\!\times_s} A \ar[r]^-{e\times id_A} \ar[dr]^{\isom}_p &
A\times A \ar[d]^m & A {{}_t\times_{id_O}} O \ar[l]_-{id_A\times e} 
\ar[dl]^q_{\isom}
\\
 & A & }
\]

Thus an internal category in \kitty{C} is defined to be a sextuple 
$C= \oset{A,O,s,t,e,m}$ where $A,\ O$ are objects (giving respectively the 
morphisms and the objects of the internal category) and $s,t,e,m$ are
morphisms satisfying the given diagrams.  These represent the category axioms, 
hence the internal categories in $\cat{Set}$ are just ordinary small
categories.  In $\cat{Grp}$, an internal category is a small category in
which both the objects and the morphisms form groups and all the structure
maps are homomorphisms.  Suppose $A$ has multiplication $\mu_A$ and $O$ has
multiplication $\mu_O$.  Then for $s$ to be a homomorphism means that the
following diagram commutes:
\[
\xymatrix{ A\times A \ar[r]^-{\mu_A} \ar[d]_{s\times s} & A \ar[d]^s \\
           O\times O \ar[r]_-{\mu_O} & O
}
\]
and similarly for the other morphisms.  We can define a multiplication on $C$
as \\ $\mu :C\times C \rightarrow C$ with $\mu = \mu_O$ on objects and $\mu =
\mu_A$ on arrows.  Then it is straightforward to check that $\mu$ is a functor 
on $C$.  Similarly, since both $O$ and $A$ are groups, these have
multiplicative identities and inverses and so we can build functors 
$\varepsilon: \mathbf{0}\rightarrow C$ (where $\mathbf{0}$ is the terminal 
object in $\cat{Cat}$, \ie\, the one-object discrete category) and 
$\iota:C\rightarrow C$ which pick out respectively an identity object and 
arrow and inverses for multiplication.  But $C$ is a small category, so it 
is an object of $\cat{Cat}$ and the functors are all morphisms 
of $\cat{Cat}$, hence from our internal category we have constructed a group 
object\footnote{The associativity of multiplication follows immediately from
  its functoriality.} in $\cat{Cat}$.  
Conversely, given a group object, $G$, in $\cat{Cat}$ we have sufficient data
to reconstruct $G$ as an internal category in $\cat{Grp}$ (essentially the
above process in reverse).  Hence internal categories in $\cat{Grp}$ are
equivalent to group objects in $\cat{Cat}$.

\subsection*{Interlude}

Suppose $A, O$ are groups with homomorphisms $\xymatrix{ A \ar@<0.5ex>[r]^s &
  O \ar@<0.5ex>[l]^e }$ such that $se = id_O$.  Then $s$ is an epimorphism
since if $a,b:O\rightarrow B$ with $as = bs$, then $ase = bse \Rightarrow
a=b$.  We call $s$ a \dfn{split epimorphism} and $e$ a \dfn{splitting} of $s$
(note that $e$ is itself a monomorphism, by much the same proof).  This
definition does not use any specific properties of groups or homomorphisms, so
it is valid in any category.

Given two groups $C, G$ with a left $G$-action on $C$, we can form a group 
$C\rtimes G = \{(c,g): c\in C, g\in G\}$ with
multiplication $(c,g)\cdot(c',g'):=(c\act{g}{c'}, gg')$.  This is like the 
direct product of $C$ and $G$, except that the component of multiplication in
$C$ is ``twisted'' by the $G$-action.  The group $C\rtimes G$ is called the
\dfn{semidirect product} of $C$ by $G$.  

Now let us return to the split epimorphism $s:A\rightarrow O$ and its splitting
$e$.  For each $a\in A$, we can write $a=ke(x)$, where $k=a(es(a))^{-1}\in 
\Ker s$ and $x = s(a)$.  Suppose $a' = k'e(x')$.  Then $aa' = ke(x)k'e(x')
=ke(x)k'(e(x))^{-1}e(x)e(x').$  We can define an action of $O$ on $\Ker s$ by
$\act{x}{k}:= e(x)k(e(x))^{-1}$ (since $e$ is a homomorphism, we can of course
rewrite $(e(x))^{-1}$ as $e(x^{-1})$).  We can then write $aa'= k
\act{x}{k'}e(xx')$.  Note that unless brackets indicate otherwise, $x$ acts
only on the symbol to its immediate right (and similarly for morphisms).  
There is a map $\phi :A\rightarrow \Ker s \rtimes O,\ \phi(ke(x))=(k,x)$.
Now $\phi(aa')=\phi(k\act{x}{k'}e(xx'))=(k\act{x}{k'},xx')=(k,x)(k',x')
=\phi(a)\phi(a')$, so $\phi$ is a homomorphism.  Also, there is an obvious
inverse $\phi^{-1}:\Ker s \rtimes O \rightarrow A,\ \phi^{-1}(k,x):=ke(x)$,
which is also a homomorphism.  Hence $\phi$ is an isomorphism and we have 
established that $A \isom \Ker s\rtimes O$.

Let $\pt :C\rightarrow G$ be a morphism of groups.  Then $\oset{C, G, \pt}$ is
a \dfn{crossed module} (of groups) if there is a (left) $G$-action on $C$ such
that the following two properties hold for all $c,d \in C,\ g\in G$:
\begin{enumerate}
  \item $\pt\act{g}{c} = g\pt cg^{-1}$ 
  \item $\act{\pt c}{d} = cdc^{-1}$
\end{enumerate}
The first property is known as \dfn{equivariance} of $\pt$ with respect to
the action (see \cite[p. 220]{KampsK:abshsh}) and the second is called the 
\dfn{Peiffer identity} (see \cite[p. 250]{PorterT:intccm})\footnote{Note that
neither of these names for the crossed module axioms is consistently applied
in the literature.  Most authors seem content just to number the axioms.}.  
For example, if $N \lhd G$ then the inclusion $N \hookrightarrow G$ gives a 
crossed module with the trivial action.  Further generic examples of crossed
modules may be found in \cite{BrownR:comhai}.  An immediate consequence of the 
crossed module axioms is that $\Ker \pt$ is abelian, while the image is a
normal subgroup of $G$.

\subsection*{Crossed Modules and Internal Categories}

We shall show that crossed modules of groups are equivalent to internal
categories in \cat{Grp}.  In other words, given any crossed module, we can
construct an internal category and vice versa.

Let $\oset{C, G, \pt}$ be a crossed module.  Since $G$ acts on $C$, we can
form the semidirect product $C \rtimes G$ as defined above and define maps
$s,t:C\rtimes G \rightarrow G$ and $e:G\rightarrow C\rtimes G$ by $s(c,g):=g,\
t(c,g):= \pt cg$ and $e(g):=(1_C,g)$.  Then $s$ is clearly a homomorphism.
Also,
\[\begin{array}{rclclcl}
 t((c,g)\cdot(d,h)) & = & t(c\act{g}{d},gh) & = & \pt(c\act{g}{d})gh & = & 
\pt c g\pt d g^{-1}gh \\
   & = & \pt cg\pt d h & = & t(c,g)t(d,h) &  &   \end{array}
\]
and $e(gh) = (1_C,gh) = (1_C,g)\cdot(1_C,h) = e(g)\cdot e(h)$
so $t$ and $e$ are also homomorphisms (for $e$ we implicitly used
the fact that the action of $G$ on $C$ determines a map  
$G \rightarrow \textit{Aut}(C)$ and hence $\act{g}{1_C}=1_C,\ \forall 
g\in G$).  Furthermore, for all $g\in G,\ se(g)=s(1_C,g)=g=id_G(g)$ and
$te(g) = t(1_C,g)=\pt(1_C)g = g$ so $se = te = id_G$, i.e. $s$ and $t$ are 
split epimorphisms with common splitting $e$.

We have thus constructed an internal (reflexive) digraph ${
  \xymatrix{ C\rtimes G \ar@<1ex>[r] \ar[r] & G \ar@/^/[l] } }$.  The
elements can be pictured as follows:
\[ 
\xymatrix{ g \ar[r]^-{(c,g)} & \pt cg }
\]

Vertices correspond to elements of $G$ and edges to elements of $C\rtimes G$.
The source and target vertices of a given edge $(c,g)$ are given by 
$s(c,g)$ and $t(c,g)$ respectively and $e(g)$ gives a loop on vertex $g$.
There is an obvious composition of edges:
\[
\xymatrix{ g \ar[r]^-{(c,g)} \ar@/_1.5pc/[rr]_-{(c'c,g)} & \pt cg \ar[r]^-{(c',\pt
    cg)} & \pt c'\pt cg .}
\]
Thus we define $(c',\pt cg)\comp (c,g):= (c'c, g)$ (note that $s(c', \pt cg) =
\pt cg = t(c,g)$).  Since $\pt$ is a homomorphism, we have $\pt(c'c)=\pt c' \pt c$
as required.

To get an internal category we now just require that $\comp$ be a
homomorphism.  In other words, we need 
\begin{equation} 
((c',\pt cg)\cdot(d',\pt dh))\comp((c,g)\cdot(d,h))=((c',\pt
cg)\comp(c,g))\cdot ((d',\pt dh)\comp(d,h)) \end{equation}
This is the familiar interchange law (see page \pageref{inter}).  Evaluating 
the two sides separately, we get:

\[\begin{array}{rclcl}
  \textrm{LHS} & = & (c'\act{\pt cg}d',\pt cg\pt dh)\comp(c\act{g}{d},gh) & = &
(c'\act{\pt cg}{d'}c\act{g}{d},gh) \\
  & = & (c'c\act{g}{d'}c^{-1}c\act{g}{d},gh) & = & (c'c\act{g}{d'}\act{g}{d},gh) \\
 & = & (c'c\act{g}{(d'd)},gh) & &
\end{array} 
\]

and 
\[\textrm{RHS} = (c'c,g)\cdot(d'd,h) = (c'c\act{g}{(d'd)},gh) \]
whence equality.  Therefore $\comp$ is indeed a homomorphism and so we have
constructed an internal category in \cat{Grp}.

Conversely, suppose that we have an internal category $\langle 
\small{\xymatrix{ A \ar@<1ex>[r]^s \ar@<0.3ex>[r]_t & O \ar@<1ex>@/^1ex/[l]^e }},\comp \rangle$.  We have seen that $A\isom \Ker s\rtimes O$ with
$O$ acting on $\Ker s$ by $\act{x}{k}:=e(x)ke(x^{-1})$ for $x\in O,
k\in \Ker s$.  Objects of the category are the elements of $O$ while morphisms
are of the form $(k,x)$ with $k\in \Ker s, x\in O$.  The maps $s,t$ give
respectively the source and target objects of each morphism (note that the
$x$ in $(k,x)$ is effectively a label to show the source of the morphism),
while $e$ gives the identity arrow for each object.

Define $\pt :\Ker s\rightarrow O$ to be the restriction of $t$ to $\Ker s$, \ie
$\pt = t\vert_{\small\Ker s}$.  Then $\pt$ is automatically a homomorphism, since
$t$ is one.  For any $k\in\Ker s,\, x\in O$, we have $\pt(\act{x}{k}) = 
t(\act{x}{k}) = t(e(x)ke(x^{-1}))$ by definition.  But $t$ is a homomorphism 
so this is the same as $te(x)tkte(x^{-1}) = xt(k)x^{-1}$ (since $te=id_O$) $=
x\pt kx^{-1}$.  Hence $\pt$ is equivariant with respect to the action.  It
remains to verify the Peiffer identity.  We know that composition is a
morphism and hence with the multiplication in $\Ker s \rtimes O$ it satisfies
the interchange law, \ie we have:
\[
((k', \pt (k)x)\cdot(l',\pt (l)y))\comp((k,x)\cdot(l,y))=
((k',\pt (k)x)\comp(k,x))\cdot ((l',\pt (l)y)\comp(l,y)).
\]
Evaluating the two sides of this equation gives:
\[
\begin{array}{rcl}
  \textrm{LHS} & = & (k'\act{\pt kx}{l'}, \pt kx\pt ly)\comp(k\act{x}{l},xy) \\
               & = & (k'\act{\pt kx}{l'}k\act{x}{l},xy)
\end{array}\]
(this composition is defined since $\pt (k\act{x}{l})xy = 
\pt k\pt\act{x}{l}xy = \pt k x\pt lx^{-1}xy = \pt kx\pt ly$ by equivariance) and
\[\textrm{RHS} = (k'k,x)\cdot(l'l,y) = (k'k\act{x}{(l'l)},xy). \]
Since the two sides are equal, we know that their first components must
be equal.  So we have
\[
\begin{array}{rcl}
 k'\act{(\pt kx)}{l'}k\act{x}{l} & = & k'\act{\pt k}{(\act{x}{l'})}k\act{x}{l}
 \\
      & = & k'k\act{x}{(l'l)} \\
      & = & k'k\act{x}{l'}\act{x}{l} \\
      & = & k'k\act{x}{l'}k^{-1}k\act{x}{l}.
\end{array} \]
Cancelling on both sides and writing $m = \act{x}{l'}\in \Ker s$, we get
$\act{\pt k}{m} = kmk^{-1}$, which is the Peiffer identity as required.  Hence
$\oset{\Ker s, O, \pt}$ is a crossed module with the action arising from the
semidirect product $\Ker s \rtimes O$.  

We have shown that crossed modules of groups are equivalent to internal
categories in $\cat{Grp}$.  We saw earlier that these are in turn equivalent
to group objects in $\cat{Cat}$ and hence have arrived at a result proved 
by Brown and Spencer \cite{BrownR:grocmf} in the 1970s, namely that crossed 
modules of groups are equivalent to group objects in $\cat{Cat}$ (they did
not go via internal categories of groups but went directly between the two
using a functor).  In fact, they went slightly further and showed that if a
category is equipped with a group structure it must in fact be a groupoid and
hence the group objects in $\cat{Cat}$ are the same as the group objects in 
$\cat{Grpd}$, the category of groupoids.  Thus they proved the equivalence of 
crossed modules of groups and group objects in the category of groupoids. 

Since doing this work I have gained access to the new edition of MacLane's
seminal volume \cite{MacLaneS:catwm}, which now includes some material on
internal categories as well as on group objects and other internal algebraic
structures.  The reader may also find it useful to consult the following
sources, which treat various aspects of this area
\cite{AlpM:enucog,AznarGarciaE:cohnac,EllisG:cromhd,KampsK:cattg,LavendhommeR:cohnsa,LodayJ:spafmn,NorrieK:cromag,PorterT:cromic,PorterT:extcma}.

\bibliographystyle{magnus}        
\bibliography{cgrefs}

\end{document}